\documentclass[a4paper,12pt]{article}
\usepackage[cp1251]{inputenc}
\usepackage[english]{babel}
\usepackage[tbtags]{amsmath}
\usepackage{amsfonts,amssymb,mathrsfs,amscd}
\usepackage{wrapfig}
\usepackage[dvips]{graphicx}
\newtheorem{theorem}{Theorem}[section]
\newtheorem{lemma}{Lemma}[section]
\newtheorem{definition}{Definition}[section]

\date{}
 \numberwithin{equation}{section}

\begin{document}

\large
 \centerline{\bf Trees in Wavelet analysis on  Vilenkin groups} \centerline{\bf S.\,F.~Lukomskii}

\noindent
 N.G.\ Chernyshevskii Saratov State University\\
 LukomskiiSF@info.sgu.ru

%\markright{Multiresolution analysis on zero-dimensional Abelian
%groups}

\footnotetext[0]{This research was carried out with the financial
support of  the Russian Foundation for Basic Research (grant
no.~13-01-00102).}

\begin{abstract}
 We consider a class of
 $(1,M)$-elementary step functions on the $p$-adic Vilenkin group.
 We prove that $(1,M)$-elementary step
 function generates a MRA on $p$-adic Vilenkin group iff it is
 generated by a rooted tree on the set of vertices
 $\{0,1,\dots p-1\}$ with 0 as a root.
Bibliography: 14 titles.
\end{abstract}
\noindent
 keywords: zero-dimensional group, Vilenkin group, multiresolution
analysis, wavelet bases, tree.

%---------------------------------------------

\section{Introduction}\label{s1}

  In articles \cite{1}-\cite{4} first
 examples of orthogonal wavelets on the dyadic Cantor  group ($p=2$) are
 constructed and their properties are studied.
  Yu.Farkov \cite{5}-\cite{7} found necessary and
 sufficient conditions for a refinable function to generate an orthogonal
 MRA in the $L_2(\mathfrak G)$
 -spaces on the $p$-adic Vilenkin group $\mathfrak G$. These conditions use the
 Strang-Fix and the modified Cohen  properties.

In \cite{7}  this
 construction  is  given in a concrete fashion for p = 3.
 In \cite{8}, some algorithms for constructing orthogonal and biorthogonal compactly
 supported wavelets on Vilenkin groups are proposed. In \cite{5}-\cite{8} two
 types of orthogonal wavelet examples are constructed: step
 functions and sums of Vilenkin  series.

 Khrennikov, Shelkovich, and Skopina
\cite{10},\cite{11} introduced the concept of a~$p$-adic MRA with
orthogonal refinable function, and described a~general pattern for
their construction. This method was developed for an orthogonal
refinable function $\varphi$ with condition ${\rm
supp}\widehat{\varphi}\subset B_0(0)$, where $B_0(0) = \{x :
|x|_p\le 1 \}$ is the unit ball in the field $\mathbb Q_p$.
Similar results were obtained for  arbitrary zero-dimensional
group \cite{13}. The condition  ${\rm
supp}\widehat{\varphi}\subset B_0(0)$ is very important. S.
Albeverio, S. Evdokimov, M. Skopina \cite{12} proved that if a
refinable step function $\varphi$ generates an orthogonal $p$-adic
MRA, then ${\rm supp}\hat \varphi(\chi)\subset B_0(0)$.

On the other hand on Vilenkin groups Yu.A.Farkov constructs
examples of step refinable functions $\varphi$ generating an
orthogonal MRA with ${\rm supp}\widehat{\varphi}\subset G_1^\bot$.
In the author's work \cite{14} a necessary condition for a support
of orthogonal refinable step function are found: if step refinable
 $(1,N)$-elementary functions $\varphi$ generated an orthogonal MRA
 on $p$-adic Vilenkin group, then ${\rm supp}\widehat{\varphi}\subset
 G_{p-2}^\bot$.

 In this work we study a structure of the set
 ${\rm supp}\widehat{\varphi}$. We prove that $(1,N)$-elementary function $\varphi$
 generates an orthogonal MRA
 on $p$-adic Vilenkin group iff the function $\varphi$  is
 generated by means of some tree. For any tree we give an
 algorithm for constructing corresponding refinable
 function and orthogonal wavelets.

The paper is organized as follows. We  consider $p$-adic Vilenkin
group $\mathfrak G$ as a zero-dimensional group $(G,\dot + )$ with
condition  $pg_n$=0. Therefore, in section 2, we recall some
concepts and facts from the theory of  zero-dimensional group. We
will systematically use the notation and the results from
\cite{13},\cite{14}.

In section 3 and the following sections we consider MRA on
$p$-adic Vilenkin group $\mathfrak G$. In section 3 we study
refinable step-functions which generate the orthogonal MRA. We
define a class of $(N,M)$-elementary set and prove that the shifts
system $\varphi (x-\dot h)_{h\in H_0}$ is orthonormal if ${\rm
supp}\widehat{\varphi}$ is $(N,M)$-elementary set.

In section 4 we introduce such  concepts as "a set generated by a
tree" and "a refinable step function generated by a tree" and
prove, that any rooted tree generates a refinable step function
that generate an orthogonal MRA on Vilenkin group.

In section 5 we give an algorithm for constructing orthogonal
wavelets according to the tree.

\section{Preliminaries}\label{s2}
We will consider the Vilenkin group as a  locally compact
zero-dimensional Abelian group
 with additional condition $p_ng_n=0$. Therefore we start with some basic notions
  and facts related to analysis on zero-dimensional groups. More
  information it is possible to find in \cite{12}--\cite{14}.

Let $(G,\dot + )$~be a~locally compact zero-dimensional Abelian
group with the topology generated by a~countable system of open
subgroups
$$
\cdots\supset G_{-n}\supset\cdots\supset G_{-1}\supset G_0\supset
G_1\supset\cdots\supset G_n\supset\cdots
$$
where
$$
\bigcup_{n=-\infty}^{+\infty}G_n= G,\quad
 \quad \bigcap_{n=-\infty}^{+\infty}G_n=\{0\},
$$
 $p_n$ be an order of quotient group $G_n/G_{n+1}$.
 We will always
assume that all~$p_n$ are prime numbers. We will name such chain
as \it basic chain. \rm  In this case, a~base of the topology is
formed by all possible cosets~$G_n\dot + g$, $g\in G$.

We further define the numbers $(\mathfrak
m_n)_{n=-\infty}^{+\infty}$ as follows:
$$
\mathfrak m_0=1,\qquad \mathfrak m_{n+1}=\mathfrak m_n\cdot p_n.
$$
Let~$\mu$ be a Haar measure on~$G$, we know that $\mu
G_n=\frac{1}{\mathfrak m_n}$. Further, let
$\smash[b]{\displaystyle\int_{G}f(x)\,d\mu(x)}$~be the absolutely
convergent integral of the measure~$\mu$.

Given an $n\in\mathbb Z$, take an element $g_n\in G_n\setminus
G_{n+1}$ and fix~it. Then any $x\in G$ has a~unique representation
of the form
\begin{equation}
\label{eq1.1} x=\sum_{n=-\infty}^{+\infty}a_ng_n, \qquad
a_n=\overline{0,p_n-1}.
\end{equation}
The sum~\eqref{eq1.1} contain finite number of terms with negative
subscripts, that~is,
\begin{equation}
\label{eq1.2} x=\sum_{n=m}^{+\infty}a_ng_n, \qquad
a_n=\overline{0,p_n-1}, \quad a_m\ne 0.
\end{equation}
We will name system $(g_n)_{n\in \mathbb Z}$ as {\it a basic
system}.

Classical examples of zero-dimensional groups are Vilenkin groups
and groups of $p$-adic numbers (see~\cite[Ch.~1, \S\,2]{12}).
 A direct sum of cyclic groups $Z(p_k)$ of order~$p_k$,
 $k\in\mathbb Z$, is called a~\textit{Vilenkin group}. This means
 that the elements of a~Vilenkin group are infinite sequences
 $x=(x_k)_{k=-\infty}^{+\infty}$ such that:
 \begin{itemize}
 \item[1)] $x_k=\overline{0,p_k-1}$; \item[2)] only a~finite number
 of~$x_k$ with negative subscripts are different from zero;
 \item[3)] the group operation~$\dot+ $ is the coordinate-wise
 addition modulo~$p_k$, that~is,
 $$
 x\dot+ y=(x_k\dot+ y_k), \qquad x_k\dot+ y_k=(x_k+y_k)\ \
 \operatorname{mod}p_k.
 $$
 \end{itemize}
 A topology on such group is generated by the chain of subgroups
 $$
  G_n=\bigl\{x\in G:x=(\dots,0,0,\dots,0,x_n,x_{n+1},\dots),\
   x_\nu=\overline{0,p_\nu-1},\ \nu\ge n\bigr\}.
 $$
 The elements $g_n=(\dots,0,0,1,0,0,\dots)$ form a basic system.
 From definition of the operation $\dot+$ we have $p_ng_n=0$.
 Therefore we will name a zero-dimensional group $(G,\dot+)$ with the condition $p_ng_n=0$ as Vilenkin group.

By $X$ denote  the collection
 of the characters of a~group $(G,\dot+ )$; it is
a~group with respect to multiplication too. Also let
$G_n^\bot=\{\chi\in X:\forall\,x\in G_n\  , \chi(x)=1\}$ be the
annihilator of the group~$G_n$. Each annihilator~$ G_n^\bot$ is
a~group with respect to multiplication, and the subgroups~$
G_n^\bot$ form an~increa\-sing sequence
\begin{equation}
\label{eq1.3} \cdots\subset G_{-n}^\bot\subset\cdots\subset
G_0^\bot \subset G_1^\bot\subset\cdots\subset
G_n^\bot\subset\cdots
\end{equation}
with
$$
\bigcup_{n=-\infty}^{+\infty} G_n^\bot=X \quad\text{and} \quad
\bigcap_{n=-\infty}^{+\infty} G_n^\bot=\{1\},
$$
the quotient group $ G_{n+1}^\bot/ G_n^\bot$ having order~$p_n$.
The group of characters~$X$ is a zero-dimensional group with a
basic chain \eqref{eq1.3}. The group  may be equipped with the
topology using the chain of subgroups~\eqref{eq1.3}, the family of
the cosets $ G_n^\bot\cdot\chi$, $\chi\in X$, being taken
as~a~base of the topology. The collection of such cosets, along
with the empty set, forms the~semiring~${\mathscr X}$. Given
a~coset $ G_n^\bot\cdot\chi$, we define a~measure~$\nu$ on it by
$\nu( G_n^\bot\cdot\chi)=\nu( G_n^\bot)= \mathfrak m_n$ (so that
always $\mu( G_n)\nu( G_n^\bot)=1$). The measure~$\nu$ can be
extended onto the $\sigma$-algebra of measurable sets in the
standard way. One then forms the absolutely convergent integral
$\displaystyle\int_XF(\chi)\,d\nu(\chi)$ of this measure.

The value~$\chi(g)$ of the character~$\chi$ at an element $g\in G$
will be denoted by~$(\chi,g)$. The Fourier transform~$\widehat f$
of an~$f\in L_2( G)$  is defined~as follows
$$
\widehat f(\chi)=\int_{ G}f(x)\overline{(\chi,x)}\,d\mu(x)=
\lim_{n\to+\infty}\int_{ G_{-n}}f(x)\overline{(\chi,x)}\,d\mu(x),
$$
the limit being in the norm of $L_2(X)$. For any~$f\in L_2(G)$,
the inversion formula is valid
$$
f(x)=\int_X\widehat f(\chi)(\chi,x)\,d\nu(\chi)
=\lim_{n\to+\infty}\int_{ G_n^\bot}\widehat
f(\chi)(\chi,x)\,d\nu(\chi);
$$
here the limit also signifies the convergence in the norm of~$L_2(
G)$. If $f,g\in L_2( G)$ then the Plancherel formula is valid
$$
\int_{ G}f(x)\overline{g(x)}\,d\mu(x)= \int_X\widehat
f(\chi)\overline{\widehat g(\chi)}\,d\nu(\chi).
$$
\goodbreak

Provided with this topology, the group of characters~$X$ is
a~zero-dimensional locally compact group; there is, however,
a~dual situation: every element $x\in G$ is a~character of the
group~$X$, and~$ G_n$ is the annihilator of the group~$ G_n^\bot$.
The union of disjoint sets $E_j$ we will denote by $\bigsqcup
E_j$.

 For any $n\in \mathbb Z$ we choose a character $r_n\in  G_{n+1}^{\bot}\backslash G_n^{\bot}$
 and fixed it. $(r)_{n\in \mathbb Z}$ is called a Rademacher system. Let us denote
  $$
    H_0=\{h\in G: h=a_{-1}g_{-1}\dot+a_{-2}g_{-2}\dot+\dots \dot+ a_{-s}g_{-s}, s\in \mathbb
    N,\ a_j=\overline{0,p-1}\},
  $$
  $$
    H_0^{(s)}=\{h\in G: h=a_{-1}g_{-1}\dot+a_{-2}g_{-2}\dot+\dots \dot+
    a_{-s}g_{-s},\ a_j=\overline{0,p-1}
    \},s\in \mathbb N.
  $$
  The set $H_0$ is an analog of the set $\mathbb N_0=\mathbb N\bigsqcup \{0\}$.

 If in the zero-dimensional group $G\ p_n=p$ for any $n\in\mathbb Z$
 then we can define the mapping ${\cal
 A}\colon G\to G$ by
 ${\cal A}x:=\sum_{n=-\infty}^{+\infty}a_ng_{n-1}$, where
 $x=\sum_{n=-\infty}^{+\infty}a_ng_n\in G$.  The mapping~${\cal A}$ is called
 a  dilation operator if~${\cal A}(x\dot+ y)={\cal A}x\dot + {\cal A}y$ for all
 $x,y\in G$. By definition, put $(\chi {\cal A},x)=(\chi, {\cal
 A}x)$.

  \begin{lemma}[\cite{14}]
 For any zero-dimensional group\\
 1) $\int\limits_{G_0^\bot}(\chi,x)\,d\nu(\chi)={\bf 1}_{G_0}(x)$,
 2) $\int\limits_{G_0}(\chi,x)\,d\mu(x)={\bf 1}_{G_0^\bot}(\chi)$.\\
 \end{lemma}

 \begin{lemma}[\cite{14}]
  If $p_n=p$ for any $n\in \mathbb Z$ and the mapping ${\cal A}$ is additive then \\
  1) $\int\limits_{G_n^\bot}(\chi,x)\,d\nu(\chi)=p^n{\bf
  1}_{G_n}(x)$,\\
  2) $\int\limits_{G_n}(\chi,x)\,d\mu(x)=\frac{1}{p^n}{\bf
  1}_{G_n^\bot}(\chi)$.
  \end{lemma}
\begin{lemma}[\cite{14}]
Let  $\chi_{n,s}=r_n^{\alpha_n}r_{n+1}^{\alpha_{n+1}}\dots
r_{n+s}^{\alpha_{n+s}}$ be a character does not belong to
$G_n^\bot$. Then
$$
\int\limits_{G_n^\bot\chi_{n,s}}(\chi,x)\,d\nu(\chi)=p^n(\chi_{n,s},x){\bf
1}_{G_n}(x).
$$
\end{lemma}
\begin{lemma}[\cite{14}]
Let
$h_{n,s}=a_{n-1}g_{n-1}\dot+a_{n-2}g_{n-2}\dot+\dots\dot+a_{n-s}g_{n-s}\notin
G_n$. Then
$$
\int\limits_{G_n\dot+h_{n,s}}(\chi,x)\,d\mu(x)=\frac{1}{p^n}(\chi,h_{n,s}){\bf
1}_{G_n^\bot}(\chi).
$$
\end{lemma}
 \begin{definition}[\cite{14}]
Let $M,N\in\mathbb N$.
 Denote by  ${\mathfrak D}_M(G_{-N})$ the set of step-functions
 $f\in L_2(G)$ such that 1)${\rm supp}\,f\subset G_{-N}$, and 2)
 $f$ is constant on cosets $G_M\dot+g$. Similarly is defined ${\mathfrak
 D}_{-N}(G_{M}^\bot)$.
 \end{definition}
\begin{lemma}[\cite{14}]
 Let $M,N\in\mathbb N$. $f\in \mathfrak D_M(G_{-N})$ if and only if $\hat f\in \mathfrak
 D_{-N}(G_M^\bot)$.
\end{lemma}
\section{MRA and refinable function on Vilenkin groups}\label{s3}

 In what follows we will consider groups $G$ for which
$p_n=p$
  and $pg_n=0$ for any $n\in \mathbb Z$. We know that it is a
  Vilenkin group. We will denote a Vilenkin group as $\mathfrak G$.

  In this group  we can choose Rademacher functions
  in various ways.
  We define Rademacher functions by the equation
  $$
  \left(r_n,\sum_{k\in\mathbb Z}a_kg_k\right)=\exp\left(\frac{2\pi i}{p}a_n\right).
  $$
  In this case
  $$
  (r_n,g_k)=\exp\left(\frac{2\pi i}{p}\delta_{nk}\right).
  $$
 Our main objective is to
 find a simple  algorithm  to get a refinable step-function that generates an orthogonal MRA on Vilenkin group.
 \begin{definition}
  A family of closed subspaces $V_n$, $n\in\mathbb Z$,
 is said to be a~multi\-resolution analysis  of~$L_2(\mathfrak G)$
 if the following axioms are satisfied:
 \begin{itemize}
 \item[A1)] $V_n\subset V_{n+1}$;
 \item[A2)] ${\vrule width0pt
 depth0pt height11pt} \overline{\bigcup_{n\in\mathbb
 Z}V_n}=L_2(\mathfrak G)$ and $\bigcap_{n\in\mathbb Z}V_n=\{0\}$;
 \item[A3)] $f(x)\in V_n$  $\Longleftrightarrow$ \ $f({\cal A} x)\in V_{n+1}$ (${\cal A}$~is a~dilation
 operator);
 \item[A4)] $f(x)\in V_0$ \ $\Longrightarrow$ \
 $f(x\dot - h)\in V_0$ for all $h\in H_0$; ($H_0$ is analog of $\mathbb
 Z$).
  \item[A5)] there exists
 a~function $\varphi\in L_2(\mathfrak G)$ such that the system
 $(\varphi(x\dot - h))_{h\in H_0}$ is an orthonormal basis
 for~$V_0$.
\end{itemize}

 A function~$\varphi$ occurring in axiom~A5 is called
a~\textit{scaling function}.
\end{definition}

 Next we will follow the conventional approach. Let
 $\varphi(x)\,{\in}\, L_2(\mathfrak G)$, and suppose that
 $(\varphi(x\dot -\nobreak h))_{h\in H_0}$ is an~orthonormal
 system in~$L_2(\mathfrak G)$. With the function~$\varphi$ and the
 dilation operator~${\cal A}$, we define the linear subspaces
 $L_n=(\varphi({\cal A}^nx\dot - h))_{h\in H_0}$ and
 closed subspaces $V_n=\overline{L_n}$. It is evident that the functions
  $p^{\frac{n}{2}}\varphi({\cal A^n}x \dot-h)_{h\in H_0}$ form
 an orthonormal basis for $V_n$, $n\in \mathbb Z$.   If subspaces $V_n$ form
 a~MRA, then the function~$\varphi$ is said to \textit{generate}
 an~MRA in~$L_2(\mathfrak G)$. If a function $\varphi$ generates an MRA, then we obtain from the axiom A1
\begin{equation}
  \label{eq3.1}
  \varphi(x)=\sum_{h\in H_0}\beta_h\varphi({\cal
  A}x\dot-h)\;\;\left(\sum|\beta_h|^2<+\infty\right).
 \end{equation}
  Therefore we will look up a~function
 $\varphi\in L_2(\mathfrak G)$, which generates an~MRA
 in~$L_2(\mathfrak G)$, as a~solution of the refinement
 equation (\ref{eq3.1}), A solution of refinement equation (\ref{eq3.1}) is called a {\it refinable function}.
 \begin{lemma}[\cite{14}]
Let $\varphi \in \mathfrak D_M(\mathfrak G_{-N})$ be a solution of
(\ref{eq3.1}). Then
\begin{equation} \label{eq3.2}
\varphi(x)=\sum_{h\in H_0^{(N+1)}}\beta_h\varphi({\cal A}x\dot-h)
 \end{equation}
\end{lemma}
 The refinement equation (\ref{eq3.2}) may be written in the form
 \begin{equation}                                      \label{eq3.3}
 \hat\varphi(\chi)=m_0(\chi)\hat\varphi(\chi{\cal
  A}^{-1}),
 \end{equation}
  where

 \begin{equation}                                      \label{eq3.4}
 m_0(\chi)=\frac{1}{p}\sum_{h\in
 H_0^{(N+1)}}\beta_h\overline{(\chi{\cal A}^{-1},h)}
 \end{equation}
is a mask of the equation (\ref{eq3.3}).
\begin{lemma}[\cite{14}]
Let $\varphi\in\mathfrak D_M(\mathfrak G_{-N})$. Then the mask
$m_0(\chi)$ is constant on cosets $\mathfrak G_{-N}^\bot\zeta$. If
$\hat\varphi(\mathfrak G_{-N}^\bot)\neq 0$ then $m_0(\mathfrak
G_{-N}^\bot)=1$.
 \end{lemma}
  \begin{lemma}[\cite{14}]
 The mask  $m_0(\chi)$ is a periodic  function with any period
 $r_1^{\alpha_1}r_2^{\alpha_2}\dots r_s^{\alpha_s}$ $(s\in\mathbb
 N,\; \alpha_j=\overline{0,p-1},\;j=\overline{1,s})$.
 \end{lemma}
 So, if $m_0(\chi)$ is a mask of (\ref{eq3.3}) then\\
 T1) $m_0(\chi)$
 is constant on cosets  $\mathfrak G_{-N}^\bot\zeta$,\\
 T2) $m_0(\chi)$ is periodic with any period
  $r_1^{\alpha_1}r_2^{\alpha_2}\dots
 r_s^{\alpha_s}$, $\alpha_j=\overline{0,p-1}$, \\
 T3)
 $m_0(\mathfrak G_{-N}^\bot)=1$. \\
 Therefore we will assume that  $m_0$
 satisfies these conditions.

  \begin{theorem}[\cite{14}]
 $m_0(\chi)$ is a mask of  equation (\ref{eq3.3}) on the class $\mathfrak
 D_{-N}(\mathfrak G_M^\bot)$ if and only if
 \begin{equation} \label{eq3.5}
 m_0(\chi)m_0(\chi{\cal A}^{-1})\dots m_0(\chi{\cal A}^{-M-N})=0
 \end{equation}
 on $\mathfrak G_{M+1}^\bot\setminus \mathfrak G_M^\bot$.
 If, in addition, the system $\varphi(x\dot-h)_{h\in H_0}$ is
 orthonormal, then $\varphi(x)$ generate an orthogonal MRA.
 \end{theorem}
So,  to find a refinable function that generates orthogonal MRA,
we need take  a function $m_0(\chi)$ that satisfies conditions T1,
T2, T3, (\ref{eq3.5}), construct the function
$$
 \hat\varphi(\chi)=\prod\limits_{k=0}^\infty m_0(\chi{\cal
 A}^{-k})\in \mathfrak D_{-N}(\mathfrak G_M^\bot)
 $$
 and check that the system $\varphi(x\dot-h)_{h\in H_0}$ is
 orthonormal.

 For any zero-dimensional group $G$ the shifts system $(\varphi(x\dot-h))_{h\in H_0}$
 is orthonormal if the condition $|\hat\varphi(\chi)|={\bf 1}_{G_0^\bot}(\chi)$ is valid \cite{14}.
 For Vilenkin group $\mathfrak G$ we can give another condition.

\begin{definition}
Let $N,M\in \mathbb N$. A set $E \subset X$ is called
$(N,M)$-elementary if $E$ is disjoint union of $p^N$ cosets
$$
\mathfrak G_{-N}^\bot\zeta_j=\mathfrak
G_{-N}^\bot\underbrace{r_{-N}^{\alpha_{-N}}r_{-N+1}^{\alpha_{-N+1}}\dots
r_{-1}^{\alpha_{-1}}}_{\xi_j}\underbrace{r_{0}^{\alpha_{0}}\dots
r_{M-1}^{\alpha_{M-1}}}_{\eta_j}=\mathfrak G_{-N}^\bot\xi_j\eta_j,
$$
$j=0,1,...,p^N-1,
j=\alpha_{-N}+\alpha_{-N+1}p+\dots+\alpha_{-1}p^{N-1}$
$(\alpha_{\nu}=\overline{0,p-1})$ such that\\
1) $\bigsqcup\limits_{j=0}^{p^N-1}\mathfrak
G_{-N}^\bot\xi_j=\mathfrak G_{0}^\bot$, $\mathfrak
G_{-N}^\bot\xi_0=\mathfrak G_{-N}^\bot$,\\
2) for any $l=\overline{0,M+N-1}$ the intersection $(\mathfrak
G_{-N+l+1}^\bot\setminus \mathfrak G_{-N+l}^\bot)\bigcap
E\ne\emptyset$.
\end{definition}
\begin{lemma}
The set $H_0\subset \mathfrak G$ is an orthonormal  system on any
$(N,M)$-elementary set $E\subset X$.
\end{lemma}
{\bf Proof.} Using the definition of $(N,M)$-elementary set we
have
$$
\int\limits_E(\chi,h)\overline{(\chi,g)}\,d\nu(x)=\sum_{j=0}^{p^N-1}\int\limits_{\mathfrak
G_{-N}^\bot\zeta_j}(\chi,h)\overline{(\chi,g)}\,d\nu(x)=
$$
$$
 =\sum\limits_{j=0}^{p^N-1}\int\limits_X{\bf 1}_{\mathfrak
G_{-N}^\bot\zeta_j}(\chi)(\chi,h)\overline{(\chi,g)}\,d\nu(x)=
$$
$$
=\sum\limits_{j=0}^{p^N-1}\int\limits_X{\bf 1}_{\mathfrak
G_{-N}^\bot\zeta_j}(\chi\eta_j)(\chi\eta_j,h)\overline{(\chi\eta_j,g)}\,d\nu(x)=
$$
 %$\zeta_j=\xi_j\eta_j$
$$
=\sum\limits_{j=0}^{p^N-1}\int\limits_X{\bf 1}_{\mathfrak
G_{-N}^\bot\xi_j}(\chi)(\chi,h)\overline{(\chi,g)}(\eta_j,h)\overline{(\eta_j,g)}\,d\nu(x).
$$
Since
$$
(\eta_j,h)=(r_0^{\alpha_0}r_1^{\alpha_1}\dots
r_{M-1}^{\alpha_{M-1}},a_{-1}g_{-1}\dot+a_{-2}g_{-2}\dot+\dots\dot+a_{-s}g_{-s})=1,
$$
$$
(\eta_j,g)=(r_0^{\alpha_0}r_1^{\alpha_1}\dots
r_{M-1}^{\alpha_{M-1}},b_{-1}g_{-1}\dot+b_{-2}g_{-2}\dot+\dots\dot+b_{-s}g_{-s})=1,
$$
then
$$
\int\limits_E(\chi,h)\overline{(\chi,g)}\,d\nu(x)=\sum\limits_{j=0}^{p^N-1}\int\limits_{\mathfrak
G_{-N}^\bot\xi_j}(\chi,h)\overline{(\chi,g)}\,d\nu(x)=\int\limits_{\mathfrak
G_{0}^\bot}(\chi,h)\overline{(\chi,g)}\,d\nu(x)=
$$
$=\delta_{h,g}$. {$\square$}
\begin{theorem}
Let $(\mathfrak G,\dot+)$ be an $p$-adic Vilenkin group,
$E\subset\mathfrak G_M^\bot$ an $(N,M)$-elementary set. If\
 $|\hat\varphi(\chi)|={\bf 1}_E(\chi)$ on $X$ then the system of
shifts $(\varphi(x\dot-h))_{h\in H_0}$ is an orthonormal system on
$\mathfrak G$.
\end{theorem}
{\bf Proof.} Let $\tilde H_0\subset H_0$ be an finite set. Using
the Plansherel equation we have
$$
\int\limits_{\mathfrak
G}\varphi(x\dot-g)\overline{\varphi(x\dot-g)}\,d\mu(x)=
\int\limits_X|\hat\varphi(\chi)|^2\overline{(\chi,g)}(\chi,h)d\nu(\chi)=
\int\limits_E(\chi,h)\overline{(\chi,g)}d\nu(\chi)=
$$
$$
=\sum_{j=0}^{p^N-1}\int\limits_{\mathfrak
G_{-N}^\bot\zeta_j}(\chi,h)\overline{(\chi,g)}\,d\nu(\chi).
$$
Transform the inner integral
$$
\int\limits_{\mathfrak
G_{-N}^\bot\zeta_j}(\chi,h)\overline{(\chi,g)}\,d\nu(\chi)=\int\limits_X{\bf
1}_{\mathfrak
G_{-N}^\bot\zeta_j}(\chi)(\chi,h)\overline{(\chi,g)}\,d\nu(\chi)=
$$
$$
=\int\limits_X{\bf 1}_{\mathfrak
G_{-N}^\bot\zeta_j}(\chi\eta_j)(\chi\eta_j,h\dot-g)\,d\nu(\chi)=\int\limits_X{\bf
1}_{\mathfrak
G_{-N}^\bot\xi_j}(\chi)(\chi\eta_j,h\dot-g)\,d\nu(\chi)=
$$
$$
=\int\limits_{\mathfrak
G_{-N}^\bot\xi_j}(\chi\eta_j,h\dot-g)\,d\nu(\chi).
$$
Repeating the arguments of lemma 3.4 we obtain
$$
\int\limits_{\mathfrak
G}\varphi(x\dot-h)\overline{\varphi(x\dot-g)}\,d\mu(x)=\delta_{h,g}.\;\;\square
$$
\begin{theorem}[\cite{14}]
 Let  $\varphi(x)\in {\mathfrak
 D}_M(\mathfrak G_{-N})$. A shifts system
 $(\varphi(x\dot-h))_{h\in H_0}$ will be orthonormal if and only
if for any
$\alpha_{-N},\alpha_{-N+1},\dots,\alpha_{-1}=\overline{(0,p-1)}$
 \begin{equation} \label{eq36}
 \sum_{\alpha_{0},\alpha_1,\dots,\alpha_{M-1}=0}^{p-1}|\hat\varphi(\mathfrak G_{-N}^\bot
 r_{-N}^{\alpha_{-N}}\dots r_0^{\alpha_0}\dots
 r_{M-1}^{\alpha_{M-1}})|^2=1.
 \end{equation}
 \end{theorem}
 \begin{lemma}[\cite{14}]
Let $\hat\varphi\in \mathfrak D_{-N}(\mathfrak G_M^\bot)$ be a
solution of the refinement equation
$$
\hat\varphi(\chi)=m_0(\chi)\hat\varphi(\chi{\cal A}^{-1})
$$
and $(\varphi(x\dot-h))_{h\in H_0}$ be an orthonormal system.\\
 Then for any
$\alpha_{-N},\alpha_{-N+1},\dots,\alpha_{-1}=\overline{0,p-1}$
 \begin{equation} \label{eq37}
\sum_{\alpha_0=0}^{p-1}|m_0(\mathfrak G_{-N}^\bot
r_{-N}^{\alpha_{-N}}r_{-N+1}^{\alpha_{-N+1}}\dots
r_{-1}^{\alpha_{-1}}r_{0}^{\alpha_{0}})|^2=1.
 \end{equation}
\end{lemma}
\section{Trees and refinable functions}\label{s4}
In this section we reduce the problem of construction of step
refinable function to construction of some tree.

 We will consider some
special class of refinable functions $\varphi(\chi)$ for which
$|\hat\varphi(\chi)|$ is a characteristic function of a set.
Define this class.
\begin{definition}
A mask $m_0(\chi)$ is called $N$-elementary $(N\in\mathbb N_0)$ if
$m_0(\chi)$ is constant on cosets $\mathfrak G_{-N}^\bot\chi$, its
modulus $m_0(\chi)$ has two values only: 0 and 1, and
$m_0(\mathfrak G_{-N}^\bot)=1$. The refinable function
$\varphi(x)$ with Fourier transform
$$
\hat\varphi(\chi)=\prod\limits_{n=0}^\infty m_0(\chi{\cal A}^{-n})
$$
is called $N$-elementary too. $N$-elementary function $\varphi$
 is called $(N,M)$-elementary if  $\hat\varphi(\chi)\in \mathfrak D_{-N}(\mathfrak
 G_{M}^\bot)$. In this case  we will
 call the Fourier transform $\hat\varphi(\chi)$  $(N,M)$-elementary, also.
\end{definition}

\begin{definition}
Let $\tilde E=\bigsqcup\limits_{\alpha_{-1},\alpha_0}\mathfrak
G_{-1}^\bot r_1^{\alpha_{-1}}r_0^{\alpha_{0}}\subset \mathfrak
G_1^\bot$ be an $(1,1)$-elementary set. We say that the set
$\tilde E_X$ is a periodic extension of $\tilde E$ if
$$
\tilde E_X=\bigcup
\limits_{s=1}^\infty\bigsqcup\limits_{\alpha_1,\dots,\alpha_s=0}^{p-1}\tilde
E r_1^{\alpha_1}r_2^{\alpha_2}\dots r_s^{\alpha_s}.
$$
We say that the set $\tilde E$ generates an $(1,M)$ elementary set
$E$, if $\bigcap\limits_{n=0}^\infty \tilde E_X{\cal A}^n=E$.
\end{definition}
Since $\tilde E_X\supset \mathfrak G_{-N}^\bot$ then
$\bigcap\limits_{n=0}^{M+1} \tilde E_X{\cal A}^n=E$ and
$\left(\bigcap\limits_{n=0}^{M+1} \tilde E_X{\cal
A}^n\right)\bigcap(\mathfrak G_{M+1}^\bot\setminus \mathfrak
G_{M}^\bot)=\emptyset$. The converse is also true. Since
$$
\left(\bigcap\limits_{n=0}^{M+1} \tilde E_X{\cal
A}^n\right)\bigcap(\mathfrak G_{M+1}^\bot\setminus \mathfrak
G_{M}^\bot)=\emptyset.
$$
Then we have
$$
\left(\bigcap\limits_{n=0}^{M+2} \tilde E_X{\cal
A}^n\right)\bigcap(\mathfrak G_{M+2}^\bot\setminus \mathfrak
G_{M+1}^\bot)=\tilde E_X\bigcap\left(\bigcap\limits_{n=0}^{M+1}
\tilde E_X{\cal A}^n\bigcap(\mathfrak G_{M+1}^\bot\setminus
\mathfrak G_{M}^\bot)\right){\cal A}=
$$
$$
=\tilde E_X\bigcap\emptyset=\emptyset.
$$

Let us write the set $\{0,1,\dots,p-1\}$ in the form
$$
\{0,u_1,u_2,\dots,u_q,\alpha_1,\alpha_2,\dots,\alpha_{p-q-1}\}=V,\;\;0=u_0,
$$
where $1\le q\le p-1$. We will consider the set $V$ as a set of
vertices. By
$T(0,u_1,u_2,\dots,u_q,\alpha_1,\alpha_2,\dots,\alpha_{p-q-1})=T(V)$
 we will denote a rooted tree on the set of vertices $V$, where 0 is a root,
$u_1,u_2,\dots,u_q$ are first level vertices,
$\alpha_1,\alpha_2,\dots,\alpha_{p-q-1}$ are remaining
vertices.\\
For example for $p=7, q=2, u_1=3, u_2=5$ we have trees

\begin{picture}(70,50)
  \put(30,10){\circle{6}}
  \put(28,8){$0$}
   \put(28,12){\vector(-1,1){6}}
   \put(32,12){\vector(1,1){6}}
  \put(20,20){\circle{6}}
  \put(18,18){$3$}
   \put(18,22){\vector(-1,1){6}}
   \put(22,22){\vector(1,1){6}}
  \put(20,20){\circle{6}}
  \put(40,20){\circle{6}}
  \put(38,18){$5$}
  \put(42,22){\vector(1,1){6}}
  \put(10,30){\circle{6}}
  \put(8,28){$1$}
  \put(30,30){\circle{6}}
  \put(28,28){$2$}
  \put(50,30){\circle{6}}
  \put(48,28){$4$}
  \put(52,32){\vector(1,1){6}}
  \put(60,40){\circle{6}}
  \put(58,38){$6$}
\end{picture}
or
\begin{picture}(70,50)
  \put(30,10){\circle{6}}
  \put(28,8){$0$}
   \put(28,12){\vector(-1,1){6}}
   \put(32,12){\vector(1,1){6}}
  \put(20,20){\circle{6}}
  \put(18,18){$3$}
   \put(18,22){\vector(-1,1){6}}
   \put(22,22){\vector(1,1){6}}
  \put(20,20){\circle{6}}
  \put(40,20){\circle{6}}
  \put(38,18){$5$}
  \put(42,22){\vector(1,1){6}}
  \put(10,30){\circle{6}}
  \put(8,28){$1$}
  \put(30,30){\circle{6}}
  \put(28,28){$2$}
  \put(50,30){\circle{6}}
  \put(48,28){$4$}
  \put(32,32){\vector(1,1){6}}
  \put(40,40){\circle{6}}
  \put(38,38){$6$}
\end{picture}

\hspace*{15mm} Figure 1 \hspace*{40mm} Figure 2\\
 and so on.

For any tree path
$P_j=(0,u_j,\alpha_{s-1},\alpha_{s-2},\dots,\alpha_{0},\alpha_{-1})$
we construct the set of cosets
\begin{equation}              \label{eq4.1}
\mathfrak G_{-1}^\bot r_{-1}^{u_j},\mathfrak G_{-1}^\bot
r_{-1}^{\alpha_{s-1}}r_{0}^{u_j},\mathfrak G_{-1}^\bot
r_{-1}^{\alpha_{s-2}}r_{0}^{\alpha_{s-1}},\dots,\mathfrak
G_{-1}^\bot r_{-1}^{\alpha_{0}}r_{0}^{\alpha_1},\mathfrak
G_{-1}^\bot r_{-1}^{\alpha_{-1}}r_{0}^{\alpha_0}.
 \end{equation}
For example for the tree on Figure 2 and the path $(0,3,2,6)$ we
have 3 cosets
$$
\mathfrak G_{-1}^\bot r_{-1}^3,\mathfrak G_{-1}^\bot r_{-1}^2
r_0^3,\mathfrak G_{-1}^\bot r_{-1}^6r_0^2,
$$
for the path $(0,3,1)$ we have two cosets
$$
\mathfrak G_{-1}^\bot r_{-1}^3,\mathfrak G_{-1}^\bot r_{-1}^1
r_0^3.
$$
We will represent the tree $T(V)$ as the tree
\begin{picture}(40,12)
  \put(20,0){\circle{6}}
  \put(18,-2){$0$}
   \put(18,2){\vector(-1,1){6}}
   \put(22,2){\vector(1,1){6}}
  %\put(20,20){\circle{6}}
  \put(10,8){$T_1$}
  \put(17,9){\circle*{1}}
  \put(21,9){\circle*{1}}
  \put(25,9){\circle*{1}}
   %\put(18,22){\vector(-1,1){6}}
   %\put(22,22){\vector(1,1){6}}
   %\put(40,20){\circle{6}}
  \put(28,8){$T_q$}
  \end{picture}
where $T_j$ are tree branches of $T(V)$ with $u_j$ as a root.
 By $E_j$ denote a union of all cosets \eqref{eq4.1} for fixed $j$
 and set
 \begin{equation}\label{eq4.2}
\tilde E=\left(\bigsqcup\limits_{j=1}^q
E_j\right)\bigsqcup\mathfrak G_{-1}^\bot.
 \end{equation}
 It is clear that $\tilde E$ is an $(1,1)$ elementary set and $\tilde E\subset \mathfrak G_{1}^\bot$.
\begin{definition}
Let $\tilde E_X$ be a periodic extension of $\tilde E$. We say
that the tree $T(V)$ generates a set $E$, if
$E=\bigcap\limits_{n=0}^\infty\tilde E_X{\cal A}^n.$
\end{definition}
\begin{lemma}
Let $T(V)$ be a rooted tree with 0 as a root. Let $E\subset X$ be
a set generated by the tree $T(V)$, $H$ a hight of $T(V)$. Then
$E$ is an $(1,H-2)$-elementary set.
\end{lemma}
{\bf Proof.} Let us denote
$$
m(\chi)={\bf 1}_{\tilde
E_X}(\chi),\;\;M(\chi)=\prod\limits_{n=0}^\infty m(\chi{\cal
A}^{-n}).
$$
First we note that $M(\chi)={\bf 1}_E(\chi)$. Indeed
$$
{\bf 1}_E(\chi)=1\Leftrightarrow\chi\in E\Leftrightarrow
\forall\,n,\;\chi{\cal A}^{-n}\in\tilde
E_X\Leftrightarrow\forall\,n,\;{\bf 1}_{\tilde E_X}(\chi{\cal
A}^{-n})=1\Leftrightarrow
$$
$$
\forall\,n,\;m(\chi{\cal
A}^{-n})=1\Leftrightarrow\prod\limits_{n=0}^\infty m(\chi{\cal
A}^{-n})=1\Leftrightarrow M(\chi)=1.
$$
It means that $M(\chi)={\bf 1}_E(\chi)$.\\
Now we will prove, that ${\bf 1}_E(\chi)=0$ for $\chi\in\mathfrak
G_{H-1}^\bot\setminus\mathfrak G_{H-2}^\bot$. Since $\tilde
E_X\supset \mathfrak G_{-1}^\bot$ it follows that ${\bf 1}_{\tilde
E_X}(\mathfrak G_{H-1}^\bot{\cal A}^{-H})={\bf 1}_{\tilde
E_X}(\mathfrak G_{-1}^\bot)=1$. Consequently
$$
\prod\limits_{n=0}^\infty{\bf 1}_{\tilde E_X}(\chi{\cal
A}^{-n})=\prod\limits_{n=0}^{H-1}{\bf 1}_{\tilde E_X}(\chi{\cal
A}^{-n})
$$
if $\chi\in \mathfrak G_{H-1}^\bot\setminus\mathfrak
G_{H-2}^\bot$. Let us denote $m(\mathfrak G_{-1}^\bot
r_{-1}^ir_0^k)=\lambda_{i+kp}$. By the definition of cosets
\eqref{eq4.1} $m(\mathfrak G_{-1}^\bot r_{-1}^ir_0^k)\ne
0\Leftrightarrow$ the pair $(k,i)$ is an edge of the tree $T(V)$.

We need prove that
$$
{\bf 1}_E(\mathfrak G_{-1}^\bot
r_{-1}^{\alpha_{-1}}r_{0}^{\alpha_{0}}\dots
r_{H-2}^{\alpha_{H-2}})=0
$$
for $\alpha_{H-2}\ne 0$. Since $\tilde E_X$ is a periodic
extension of $\tilde E$ it follows that the function $m(\chi)={\bf
1}_{\tilde E_X}(\chi)$ is periodic with any period
$r_{1}^{\alpha_{1}}r_{2}^{\alpha_{2}}\dots r_{s}^{\alpha_{s}}$,
$s\in\mathbb N$, i.e. $m(\chi
r_{1}^{\alpha_{1}}r_{2}^{\alpha_{2}}\dots
r_{s}^{\alpha_{s}})=m(\chi)$ when $\chi\in \mathfrak G_1^\bot$.
Using this fact we can write $M(\chi)$ for  $\chi\in \mathfrak
G_{\alpha_{H-1}}^\bot\setminus \mathfrak G_{\alpha_{H-2}}^\bot$ in
the form
$$
M(\mathfrak G_{-1}^\bot\zeta)=M(\mathfrak G_{-1}^\bot
r_{-1}^{\alpha_{-1}}r_{0}^{\alpha_{0}}\dots
r_{H-2}^{\alpha_{H-2}})=
$$
$$
=m(\mathfrak G_{-1}^\bot
r_{-1}^{\alpha_{-1}}r_{0}^{\alpha_{0}})m(\mathfrak G_{-1}^\bot
r_{-1}^{\alpha_{0}}r_{0}^{\alpha_{1}})\dots m(\mathfrak
G_{-1}^\bot r_{-1}^{\alpha_{H-3}}r_{0}^{\alpha_{H-2}})m(\mathfrak
G_{-1}^\bot r_{-1}^{\alpha_{H-2}})=
$$
$$
=\lambda_{\alpha_{-1}+p\alpha_0}\lambda_{\alpha_{0}+p\alpha_1}\dots\lambda_{\alpha_{H-3}+p\alpha_{H-2}}\lambda_{\alpha_{H-2}},\;\alpha_{H-2}\ne
0.
$$
If $\lambda_{\alpha_{H-2}}=0$ then $M(\mathfrak
G_{-1}^\bot\zeta)=0$. Let $\lambda_{\alpha_{H-2}}\ne 0$. It means
that $\alpha_{H-2}=u_j$ for some $j=\overline{1,q}$. If
$\lambda_{\alpha_{H-3}+p\alpha_{H-2}}=0$ then $M(\mathfrak
G_{-1}^\bot\zeta)=0$. Therefore we assume that
$\lambda_{\alpha_{H-3}+p\alpha_{H-2}}\ne 0$. It is true iff the
pair $(\alpha_{H-2},\alpha_{H-3})$ is an edge of $T(V)$. Repeating
these arguments, we obtain a path
$(0,u_j=\alpha_{H-2},\alpha_{H-3},\dots,\alpha_s)$ of the tree
$T(V)$. Since ${\rm hight}(T)= H$ it follows that $s\ge 0$.
Consequently $(\alpha_s,\alpha_{s-1})$ is not edge and
$\lambda_{\alpha_{s-1}+p\alpha_s}=0$, where $s\ge 0$. It means
that $M(\mathfrak G_{-1}^\bot\zeta)=0$.

Now we prove that $E$ is $(1,H-2)$ elementary set. Indeed, any
path
$(0,u_j=\alpha_{s-1},\alpha_{s-2},\dots,\alpha_0,\alpha_{-1})$
defines the coset $\mathfrak G_{-1}^\bot
r_{-1}^{\alpha_{-1}}r_{0}^{\alpha_{0}}\dots
r_{s-1}^{\alpha_{s-1}}\subset E$. But for any
$\alpha_{-1}=\overline{0,p-1}$ there exists unique path with
endpoint $\alpha_{-1}$ and starting point zero. It means that $E$
is $(1,H-2)$-elementary set. $\square$
\begin{theorem}
Let $M,p\in\mathbb N$, $p\ge 3$. Let $E\subset \mathfrak G_M^\bot$
be an $(1,M)$-elementary set, $\hat\varphi\in\mathfrak
D_{-1}(\mathfrak G_M^\bot)$, $|\hat\varphi(\chi)|={\bf
1}_E(\chi)$, $\hat\varphi(\chi)$ the solution of the equation
\begin{equation}              \label{eq4.3}
\hat\varphi(\chi)=m_0(\chi)\hat\varphi(\chi{\cal A}^{-1}),
\end{equation}
where $m_0(\chi)$ is a 1-elementary mask. Then there exists a
rooted tree $T(V)$ with ${\rm height}(T)=M+2$ that generates the
set $E$.
\end{theorem}
{\bf Prof.} Since the set $E$ is $(1,M)$-elementary set and
$|\hat\varphi(\chi)|={\bf 1}_E(\chi)$, it follows from theorem 3.2
that the system $(\varphi(x\dot-h))_{h\in H_0}$ is an orthonormal
system in $L_2(\mathfrak G)$. Using the theorem 3.3 we obtain that
for $\alpha_{-1}=\overline{0,p-1}$
$$
\sum_{\alpha_0,\alpha_1,\dots.\alpha_{M-1}=0}|\hat\varphi(G_{-1}^\bot
r_{-1}^{\alpha_{-1}}r_{0}^{\alpha_0}\dots
r_{M-1}^{\alpha_{M-1}})|^2=1.
$$
Since $\hat\varphi$ is a solution of refinement equation
\eqref{eq4.3} it follows from lemma 3.5 that for
$\alpha_{-1}=\overline{0,p-1}$
\begin{equation}              \label{eq4.4}
\sum_{\alpha_0=0}^{p-1}|m_0(\mathfrak G_{-1}^\bot
r_{-1}^{\alpha_{-1}}r_{0}^{\alpha_{0}})|^2=1.
\end{equation}
Let as denote $\lambda_{\alpha_{-1}+p\alpha_0}:=m_0(\mathfrak
G_{-1}^\bot r_{-1}^{\alpha_{-1}}r_{0}^{\alpha_{0}})$. Then we
write \eqref{eq4.4} in the form
\begin{equation}              \label{eq4.5}
\sum_{\alpha_0=0}^{p-1}|\lambda_{\alpha_{-1}+p\alpha_0}|^2=1.
\end{equation}
Since the mask $m_0(\chi)$ is 1-elementary it follows that
$|\lambda_{i+pj}|$ take two value only: 0 or 1.

Now we will construct the tree $T$. Let $\mathfrak U$ be a family
of cosets $\mathfrak G_{-1}^\bot\zeta\subset \mathfrak G_M^\bot$
such that $\hat\varphi(\mathfrak G_{-1}^\bot \zeta)\ne 0$ and
$\mathfrak G_{-1}^\bot\notin \mathfrak U$. We can write a coset
$\mathfrak G_{-1}^\bot\zeta\in \mathfrak U$ in the form
$$
\mathfrak G_{-1}^\bot r_{-1}^{\alpha_{-1}}r_{0}^{\alpha_{0}}\dots
r_{M-1}^{\alpha_{M-1}}.
$$
If $\mathfrak G_{-1}^\bot\zeta\in \mathfrak G_{n}^\bot\setminus
\mathfrak G_{n-1}^\bot$ $(n\le M)$ then $\mathfrak
G_{-1}^\bot\zeta=\mathfrak G_{-1}^\bot
r_{-1}^{\alpha_{-1}}r_{0}^{\alpha_{0}}\dots
r_{n-1}^{\alpha_{n-1}}$, $\alpha_{n-1}=\overline{1,p-1}$.

Let $u\in\overline{1,p-1}$. By $T_u$ we denote the set of vectors
$(u,\alpha_{n-1},\dots,\alpha_{0},\alpha_{-1})$ for which
$\mathfrak G_{-1}^\bot r_{-1}^{\alpha_{-1}}r_{0}^{\alpha_{0}}\dots
r_{n-1}^{\alpha_{n-1}}r_n^u\in \mathfrak U$. We will name the
vector\\ $(u,\alpha_{n-1},\dots,\alpha_{0},\alpha_{-1})$
 as a path
too. So $T_u$ is the set of pathes with starting point $u$, for
which $\hat\varphi(\mathfrak G_{-1}^\bot
r_{-1}^{\alpha_{-1}}r_{0}^{\alpha_{0}}\dots
r_{n-1}^{\alpha_{n-1}}r_n^u)\ne 0$. We will show that $T_u$ is a
rooted tree with $u$ as a root.

1) All vertices $\alpha_j,u$ of the path
$(u,\alpha_{n-1},\dots,\alpha_{0},\alpha_{-1})$ are pairwise
distinct. Indeed
$$
\hat\varphi(\mathfrak G_{-1}^\bot
r_{-1}^{\alpha_{-1}}r_{0}^{\alpha_{0}}\dots
r_{n-1}^{\alpha_{n-1}}r_n^u)=\lambda_{\alpha_{-1}+\alpha_0p}
\lambda_{\alpha_{0}+\alpha_1p}\dots\lambda_{\alpha_{n-1}+pu}\lambda_u\ne
0,\;u\ne0.
$$
If $\alpha_{n-1}=u$ then $|\lambda_{u+pu}|=|\lambda_{u+p0}|=1$
that
contradicts the equation \eqref{eq4.5}.\\
If $\alpha_{n-1}=0$ then $|\lambda_{0+pu}|=|\lambda_{0+p0}|=1$
that contradicts the equation \eqref{eq4.5} too. Consequently
$\alpha_{n-1}\notin  \{0,u\}$. By analogy we obtain that
$\alpha_i\notin
\{0,u,\alpha_{n-1},\dots,\alpha_{i+2},\alpha_{i+1}\}$.

2) If two patches $(u,\alpha_{n-1},\dots,\alpha_{0},\alpha_{-1})$
and $(u,\beta_{l-1},\dots,,\beta_{0},\beta_{-1})$ have the common
subpath
$(u,\alpha_{k-1},\dots,\alpha_{k-j+1},\alpha_{k-j})=(u,\beta_{l-1},\dots,\beta_{l-j+1},\beta_{l-j})$
and $\alpha_{k-j-1}\ne \beta_{l-j-1}$ then
$\{\alpha_{-1},\alpha_{0},\dots,\alpha_{k-j-1}\}\bigcap\{\beta_{-1},\beta_{0},\dots,\beta_{l-j-1}\}=\emptyset$.
Indeed, assume
$$
\{\alpha_{-1},\alpha_{0},\dots,\alpha_{k-j-1}\}\bigcap\{\beta_{-1},\beta_{0},\dots,\beta_{l-j-1}\}\ne\emptyset.
$$
Then there exists
$v\in\{\alpha_{-1},\alpha_{0},\dots,\alpha_{k-j-1}\}\bigcap\{\beta_{-1},\beta_{0},\dots,\beta_{l-j-1}\}$.\\
Assume that $v\ne \alpha_{k-j-1}$. Then $v=\alpha_\nu$, $-1\le\nu
\le k-j-2$ and $v=\beta_\mu$, $-1\le\mu\le l-j-1$. It follows that

$$
(u=\alpha_k,\dots,\alpha_{k-j},\alpha_{k-j-1},\dots,\alpha_{\nu+1},\alpha_{\nu}=
\beta_{\mu},\beta_{\mu-1},\dots,\beta_{0},\beta_{-1})\in T_u
$$

$$
(u=\beta_{l},\dots,\beta_{l-j}=\alpha_{k-j},\beta_{l-j-1},\dots,\beta_{\mu+1},\beta_\mu,\beta_{\mu-1},\dots,
\beta_0,\beta_{-1})\in T_u.
$$
So we have two different patches with the same sheet $\beta_{-1}$.
But this contradicts theorem 3.3. This means that $T_u$ has no
cycles, consequently $T_u$ is a graph with $u$ as a root.

3) By analogy we can proof that different trees $T_u$ an $T_v$ has
no common vertices. It follows that the graph
$T=(0,T_{u_1},\dots,T_{u_q})$ is a tree with 0 as a foot.

4) It is evident that this tree generates refinable function
$\hat\varphi$ with a mask $m_0$. Show that ${\rm height}(T)=M+2$.
Indeed, since $\hat\varphi\in  \mathfrak D_{-1}(\mathfrak
 G_{M}^\bot)$ it follows that there exists a coset $\mathfrak G_{-1}^\bot r_{-1}^{\alpha_{-1}}r_{0}^{\alpha_{0}}
 \dots r_{M-1}^{\alpha_{M-1}}$, $\alpha_{M-1}\ne 0$ for which $|\hat\varphi(\mathfrak G_{-1}^\bot r_{-1}^{\alpha_{-1}}r_{0}^{\alpha_{0}}
 \dots r_{M-1}^{\alpha_{M-1}})|=1$. This coset generates a path $(0,\alpha_{M-1}=u,\alpha_{M-2},\dots,
 \alpha_0,\alpha_{-1})$ of $T$. This path contain $M+2$ vertex. It
 means that ${\rm height}(T)\ge M+2$. On the other hand there isn't
 coset $\mathfrak G_{-1}^\bot\zeta\subset\mathfrak G_{M+1}^\bot\setminus \mathfrak
 G_{M}^\bot$, consequently  there isn't path with $L>M+2$. So
${\rm height}(T)=M+2$. Since ${\rm supp}\,\hat\varphi(\chi)$ is
$(1,M)$-elementary set, it follows that the set of all vertices of
the tree $T$ is the set $\{0,1,\dots,p-1\}$. The theorem is
proved.  $\square$
\begin{definition}
Let $T(V)$ be a rooted tree with 0 as a root,  $H$ a hight of
$T(V)$, $V=\{0,1,...,p-1\}$. Using cosets \eqref{eq4.1} we define
the mask $m_0(\chi)$ in the subgroup $\mathfrak G_{1}^\bot$ as
follows: $m_0(\mathfrak G_{-1}^\bot)=1, m_0(\mathfrak G_{-1}^\bot
r_{-1}^{i}r_0^{j})=\lambda_{i+pj}$, $|\lambda_{i+pj}|=1$ when
$\mathfrak G_{-1}^\bot r_{-1}^{i}r_0^{j}\subset \tilde E$,
\eqref{eq4.2}, $|\lambda_{i+pj}|=1$ when $\mathfrak G_{-1}^\bot
r_{-1}^{i}r_0^{j}\subset \mathfrak G_{1}^\bot\setminus\tilde E$.
Let us extend the mask $m_0(\chi)$ on the $X\setminus \mathfrak
G_{1}^\bot$ periodically, i,e $m_0(\chi
r_1^{\alpha_1}r_2^{\alpha_2}\dots r_s^{\alpha_s})=m_0(\chi)$.
Then we say that the tree $T(V)$ generates the mask $m_0(\chi)$.
Set $\hat \varphi(\chi)=\prod\limits_{n=0}^\infty m_0(\chi{\cal
A}^{-n})$.
It follows from lemma 4.1 that\\
 1) ${\rm
supp}\,\hat\varphi(\chi)\subset\mathfrak G_{H-2}^\bot$,\\
 2)
$\hat\varphi(\chi)$ is $(1,H-2)$ elementary function,\\
 3)
$(\varphi(x\dot-h))_{h\in H_0}$ is an orthonormal system.\\ In
this case we say that the tree $T(V)$ generates the refinable
function $\varphi(x)$.
\end{definition}

\begin{theorem}
Let $p\ge 3$ be a prime number,
$$
V=\{0,u_1,u_2,\dots,u_q,a_1,a_2,\dots,a_{p-q-1}\}
$$ a set of vertices, $T(V)$ a rooted tree, 0 the root,
$u_1,u_2,\dots,u_q$ a first level vertices. Let $H$ be are height
of $T(V)$. By $\varphi(x)$ denote the function generated by the
$T(V)$. Then $\varphi(x)$ generate an orthogonal MRA on $p$-adic
Vilenkin group.
\end{theorem}
{\bf Proof.} Since  $T(V)$ generates the the function  $\varphi$
 then 1)$\hat\varphi\in\mathfrak
D_{-1}(\mathfrak G_1^\bot)$, 2)$\hat\varphi(\chi)$ is $(1,H-2)$
elementary function, 3)$\hat\varphi(\chi)$ is a solution of
refinable equation \eqref{eq3.3}, 4)$(\varphi(x\dot-h))_{h\in
H_0}$ is an orthonormal system. From the theorem 3.1 it follows
that $\varphi(x)$ generates an orthogonal MRA. $\square$

{\bf Remark.} It is possible to give an algorithm for constructing
the refinable function $\varphi(x)$. Let $T(V)$ be a tree on the
set $\{0,1,\dots,p-1\}$. Construct a finite sequence
$(\lambda_{i+jp})_{i,j=0}^{p-1}$ as follows: $\lambda_0=1$,
$|\lambda_{i+pj}|=1$ if the pair $(j,i)$ is an edge of $T(f)$. For
any vertex $\alpha_{-1}$ we take the path
$(0=\alpha_{s+1},u_j=\alpha_s,
\alpha_{s-1},\dots,\alpha_0,\alpha_{-1})$ and suppose
$$
\hat\varphi(\mathfrak G_{-1}^\bot
r_{-1}^{\alpha_{-1}}r_{0}^{\alpha_{0}}\dots
r_{s-1}^{\alpha_{s-1}}r_{s}^{\alpha_{s}}r_{s+1}^{0})=\lambda_{\alpha_{-1}+\alpha_0p}\cdot\lambda_{\alpha_{0}+\alpha_1p}\cdot\dots
\cdot\lambda_{\alpha_{s-1}+\alpha_sp}\cdot\lambda_{\alpha_s}.
$$
Otherwise we suppose $\hat\varphi(\mathfrak G_{-1}^\bot\zeta)=0$.
Then $\varphi$ generates an orthogonal MRA on $p$-adic Vilenkin
group $\mathfrak G$.

\section{Construction of wavelet bases}\label{s5}
In \cite{6} and \cite{7} Yu.A.Farkov reduces the problem of
$p$-wavelet decomposition into a problem of matrix extension. We
will use more simple method \cite{13}.

As usual, $W_n$ stands for the orthogonal complement of $V_n$ in
$V_{n+1}$: that is $V_{n+1}=V_n\oplus W_n$ and $V_n\bot W_n$
($n\in\mathbb Z$, and $\oplus$ denotes the direct sum).\\
It is readily seen that\\
1) $f\in W_n\Leftrightarrow f({\cal A}x)\in W_{n+1}$,\\
2) $W_n\bot W_k$ for $k\ne n$,\\
3) $\oplus W_n=L_2(\mathfrak G)$, $n\in\mathbb Z$.

From theorems 4.1, 4.2 we derive an algorithm for constructing
wavelet bases.\\
{\bf Step 1.} Choose an arbitrary tree
$T(V)=T(0,u_1,\dots,u_q,\alpha_1,\dots,\alpha_{p-q-1})$ on the set
$V=\{0,1,\dots,p-1\}$. Let $H$ be a height of the tree $T(V)$.\\
{\bf Step 2.} Choose a finite sequence
$(\lambda_{i+jp})_{i,j=0}^{p-1}$ such that $\lambda_0=1,$
$|\lambda_{i+jp}|=1$ if the pair $(j,i)$ is the edge of the tree
$T(V)$, $|\lambda_{i+jp}|=0$ otherwise.\\
{\bf Step 3.} Construct the mask $m_0(\chi)$ and Fourier transform
$\hat\varphi(\chi)$ using definition 4.4. It is clear
that $E={\rm supp}(\hat\varphi(\chi))$ is $(1,H-2)$-elementary set.\\
{\bf Step 4.} Find coefficients $\beta_n$ for which
\begin{equation}              \label{eq5.1}
m_0(\chi)=\frac{1}{p}\sum_{h\in
H_0^{(2)}}\beta_h\overline{(\chi{\cal A}^{-1},h)}.
\end{equation}
To find coefficients $\beta_h$, we write this equation in the form
\begin{equation}\label{eq5.2}
m_0(\chi_k)=\frac{1}{p}\sum_{j=0}^{p^2-1}\beta_j\overline{(\chi_k,{\cal
A}^{-1}h_j)}
\end{equation}
 where
$$
\begin{array}{lll}
  h_j=a_{-1}g_{-1}\dot+a_{-2}g_{-2}, & j=a_{-1}+a_{-2}p, & a_{-1},a_{-2}=\overline{0,p-1}. \\
  \chi_k\in \mathfrak G_{-1}^\bot
r_{-1}^{\alpha_{-1}}r_{0}^{\alpha_{0}}, & k=\alpha_{-1}+\alpha_0p,
 &  \alpha_{-1},\alpha_0=\overline{0,p-1}.\\
\end{array}
$$

Since the matrix $\frac{1}{p}\overline{(\chi_k,{\cal A}^{-1}h_j)}$
of this system is unitary it follows that the system \eqref{eq5.2}
has a unique solution.\\
 {\bf Step 5.} We set $m_l(\chi)=m_0(\chi r_0^{-l})$,
 $l=\overline{1,p-1}$, $X_0=\{\chi:|m_0(\chi)|=1\}$. Clearly, $m_l(\chi)$ may be written as
 $$
m_l(\chi)=\frac{1}{p}\sum_{h\in H_0^{(2)}}\beta_h\overline{(\chi
r_0^{-l},{\cal A}^{-1}h)}=\frac{1}{p}\sum_{h\in
H_0^{(2)}}\beta_h^{(l)}\overline{(\chi,{\cal A}^{-1}h)}
 $$
where $\beta_h^{(l)}=\beta_h(r_0^l,{\cal A}^{-1}h)$. By the
construction of $\hat\varphi(\chi)$ we have $|m_l(X_0 r_0^l)|=1$,
$|m_l(X_0 r_0^\nu)|=0$ for $\nu\ne l$, $m_l(\chi)m_k(\chi)=0$ when $k\neq l$.\\
{\bf Step 6.} Define the functions
$$
\psi_l(x)=\sum_{h\in H_0^{(2)}}\beta_h^{(l)}\varphi({\cal
A}x\dot-h).
$$
\begin{theorem}
The functions $\psi_l(x\dot-h)$, where $l=\overline{1,p-1}$, $h\in
H_0$, form an orthonormal basis for $W_0$.
\end{theorem}
{\bf Proof.} a) We claim that
$(\varphi(\cdot\dot-g^{(1)}),\psi_l(\cdot\dot-g^{(2)}))=0$ for any
$g^{(1)}$,  $g^{(2)}\in H_0$. Since
$$
\hat\varphi_{\cdot\dot-h}(\chi)=\overline{(\chi,h)}\hat\varphi(\chi),\;\;\hat\varphi_{{\cal
A}\cdot\dot-g}(\chi)=\frac{1}{p}\overline{(\chi,{\cal
A}^{-1}g)}\hat\varphi(\chi{\cal A}^{-1}),
$$
it follows that
$$
(\varphi(\cdot\dot-g^{(1)}),\psi_l(\cdot\dot-g^{(2)}))=\int\limits_X\hat\varphi(\chi)\overline{\hat\varphi(\chi{\cal
A}^{-1})}\overline{(\chi,g^{(1)})}(\chi,g^{(2)})\overline{m_l(\chi)}\,d\nu(\chi)=0
$$
because ${\rm supp}\,\hat\varphi(\chi)=E$ and $m_l(E)=0$,
$l=\overline{1,p-1}$.\\
b) By analogy
$$
(\psi_k(\cdot\dot-g^{(1)}),\psi_l(\cdot\dot-g^{(2)}))=\int\limits_X|\hat\varphi(\chi{\cal
A}^{-1})|^2(\chi,g^{(2)}\dot-g^{(1)})m_k(\chi)\overline{m_l(\chi)}\,d\nu(\chi)=0
$$
when $k\ne l$.\\
c) We verify that
$(\psi_l(\cdot\dot-g^{(1)}),\psi_l(\cdot\dot-g^{(2)}))=0$,
provided that $g^{(1)},g^{(2)}\in H_0$ and $g^{(1)}\ne g^{(2)}$.
Write this scalar product in the form
$$
(\psi_l(\cdot\dot-g^{(1)}),\psi_l(\cdot\dot-g^{(2)}))=\int\limits_X|\hat\varphi(\chi{\cal
A}^{-1})|^2(\chi,g^{(2)}\dot-g^{(1)})|m_l(\chi)|^2\,d\nu(\chi)=
$$
$$
=\int\limits_{E{\cal A}\bigcap
X_0r_0^l}(\chi,g^{(2)}\dot-g^{(1)})\,d\nu(\chi).
$$
Show that $E{\cal A}\bigcap X_0r_0^l$ is an $(1,H-1)$-elementary
set. By the definition
\begin{equation}              \label{eq5.3}
E=\bigsqcup\limits_{(0,\alpha_s,\alpha_{s-1},\dots,\alpha_0,\alpha_{-1})\in
T(V)}\mathfrak G_{-1}^\bot
r_{-1}^{\alpha_{-1}}r_{0}^{\alpha_{0}}\dots
r_{s}^{\alpha_{s}}r_{s+1}^{0}\;\;(s\le H-3)
\end{equation}
where the union is taken over all paths
$(0,\alpha_s,\alpha_{s-1},\dots,\alpha_0,\alpha_{-1})$ of the tree
$T(V)$. It means that for any $\alpha_{-1}=\overline{0,p-1}$ the
union \eqref{eq5.3} contains unique coset $\mathfrak G_{-1}^\bot
r_{-1}^{\alpha_{-1}}r_{0}^{\alpha_{0}}\dots
r_{s}^{\alpha_{s}}r_{s+1}^{0}$.

Consequently
$$
E{\cal
A}=\bigsqcup\limits_{(0,\alpha_s,\alpha_{s-1},\dots,\alpha_0,\alpha_{-1})\in
T(V)}\mathfrak G_{0}^\bot
r_{0}^{\alpha_{-1}}r_{1}^{\alpha_{0}}\dots
r_{s+1}^{\alpha_{s}}r_{s+2}^{0}=
$$
$$
=\sum_{\alpha_{-2}=0}^{p-1}\bigsqcup\limits_{(0,\alpha_s,\alpha_{s-1},\dots,\alpha_0,\alpha_{-1})\in
T(V)}\mathfrak G_{-1}^\bot
r_{-1}^{\alpha_{-2}}r_{0}^{\alpha_{-1}}\dots
r_{s+1}^{\alpha_{s}}r_{s+2}^{0}.
$$
On the other hand
$$
X_0r_0^l=\bigcup\limits_{j\in\mathbb
N}\bigsqcup\limits_{(\gamma_{-1},\gamma_0)\in
T(V)}\bigsqcup\limits_{b_1,b_2,\dots,b_j=0}^{p-1}\mathfrak
G_{-1}^\bot r_{-1}^{\gamma_{-1}}r_{0}^{\gamma_{0}+l}r_1^{b_1}\dots
r_{j}^{b_j}.
$$
Therefore
$$
E{\cal A}\bigcap
X_0r_0^l=\bigsqcup\limits_{\gamma_{-1}=0}^{p-1}\bigsqcup\limits_{(0,\alpha_s,\alpha_{s-1},\dots,\alpha_{-1}=
\gamma_0+l,\gamma_{-1})\in T(V)}\mathfrak G_{-1}^\bot
r_{-1}^{\gamma_{-1}}r_{0}^{\gamma_{-1}}r_1^{\alpha_0}\dots
r_{s+1}^{\alpha_{s}}r_{s+2}^{0}.
$$
It means that $E{\cal A}\bigcap X_0r_0^l$ is $(1,H-1)$-elementary
set. By lemma 3.4 it follows that
$$
\int\limits_{E{\cal A}\bigcap
X_0r_0^l}(\chi,g^{(2)}\dot-g^{(1)})\,d\nu(\chi)=0.
$$
d) We claim that any function $f\in W_0$ can be expanded uniquely
in a series in $(\psi_l(x\dot-g))_{l=\overline{1,p-1},g\in H_0}$.
The proof of this fact may be found in [13], theorem 5.1.
$\square$\\
{\bf Step 7.} Since the subspaces $(V_j)_{j\in\mathbb Z}$ form an
MRA in $L_2(\mathfrak G)$, it follows that the functions
$$
(\psi_l({\cal A}^nx\dot-h))\;\;l=\overline{1,p-1}, n\in\mathbb Z,
h\in H_0
$$
form a complete orthogonal system in $L_2(\mathfrak G)$.

 \end{document}